\title{The cover time, the  blanket time, and the Matthews bound}

\author{J. Kahn \thanks{ Department of Mathematics,
Rutgers University, NJ 08903; jkahn@math.rutgers.edu. Supported in part by NSF.}
\and J. H.  Kim \thanks{Microsoft Research, Microsoft Corporation,
Redmond, WA 98052; jehkim@microsoft.com } \and L. Lov\'asz
\thanks{Microsoft Research, Microsoft Corporation, Redmond, WA
98052; lovasz@microsoft.com  } \and V. H. Vu
\thanks{ Microsoft Research, Microsoft Corporation, Redmond, WA 98052;
vanhavu@microsoft.com} }

\documentclass[11pt, twocolumn]{article}
\input amssymb.sty

\newtheorem{theo}{Theorem}[section]
\newtheorem{lemma}[theo]{Lemma}

\newtheorem{cor}[theo]{Corollary}      
\newtheorem{prop}[theo]{Proposition}
\newtheorem{claim}[theo]{Claim}

\newcommand{\TT}{{\cal T}}
\newcommand{\old}[1]{}
\newcommand\raf[1]{(\ref{#1})}

\def\E{{\sf E}}
\def\eps{\varepsilon}
\def\proofend{\hfill$\square$\medskip}
\def\Proof{\medskip\noindent{\bf Proof. }}
\renewcommand{\Pr}{{\sf P}}
\newcommand{\corr}[1]{}

\begin{document}
\maketitle
\begin{abstract}
We prove upper and lower bounds and give an approximation
algorithm for the cover time of the random walk on a graph. We
introduce a parameter $M$ motivated by the well known Matthews
bounds on the cover time and prove that $M/2\leq C = O(M (\ln \ln
n)^2)$. We give a deterministic polynomial time algorithm to
approximate $M$  within a factor of $2$; this then approximates
$C$ within a factor of $O((\ln \ln n)^2)$, improving previous
bound of $O(\ln n)$ of Matthews.

The blanket time $B$ was introduced by Winkler and Zuckerman: it
is the expectation of the first time when all vertices are visited
within a constant factor of number of times suggested by the
stationary distribution. Obviously $C \leq B$, and they
conjectured $B =O(C)$ and proved $B=O(C\ln n )$. Our bounds above
are also valid for the blanket time, and so it follows that $B =
O(C (\ln \ln n)^2)$.
\end{abstract}

\section{Introduction}

Given a connected graph $G$ on $n$ vertices, for a vertex $i\in
V(G)$, $C(i)$ denotes the cover time of
\corr{a random walk on $G$ ,} 
     the usual random walk on $G$,
\corr{starting from $i$. To be more precise, $C(i)$ is}
     starting from $i$; that is, $C(i)$ is
the expectation of the number of steps  a random walk starting
from $i$ takes until it covers all vertices of $G$. The quantity
$C = \max_{i \in V(G)} C(i)$ is called the cover time of $G$. (See
\cite{AF} for background.)

Although $C$ is a basic notion in the theory of random walks,
there is
\corr{no effective way}
      no effective way known
to compute this parameter, given the adjacency matrix of $G$ as
the input. The following question has been open for several years
\cite{AF}.

\smallskip

\noindent{\bf Question.} {\it Is there a deterministic algorithm
which approximates $C$ up to a constant factor in polynomial time
? }

\smallskip

The requirement that the algorithm is deterministic is crucial and
this makes the problem difficult. It is simple to provide a
randomized algorithm which approximates $C$ within a factor $(1+
\epsilon)$ for any positive constant $\epsilon$, with high
probability: just simulate the chain and take the average of the
empirical cover times.

Prior to this paper, the best approximation factor we knew of was
$\ln n$. This factor can be achieved using the following
\corr{classical}
      fundamental
result of Matthews \cite{Mat}. For any pair of vertices $i,j \in
V(G)$, $H(i,j)$ denotes the hitting time form $i$ to $j$. We set
\corr{Comment: Shouldn't we define hitting time?}

\[
h_{\max}=\max_{i,j\in V} H(i,j), \qquad h_{\min}=\min_{i,j\in V}
H(i,j),
\]
and more generally, for every set $S\subseteq V$, we let
\[
h_S=\min_{i,j\in S \atop i \not= j} H(i,j).
\]

\newcommand{\har}{{\rm har}}
 Let $ \har (n) =\sum_{i=1}^{n} 1/i . $
\begin{theo}[Matthews' theorem]\label{MATTH}
\corr{We have}
For any $G$,
\[
h_{\min}\har(n) \le C  \le h_{\max}\har(n).
\]
More generally, for any subset $S \subset V(G)$ with $|S|\ge 2$,
\begin{equation}\label{MATTHL}
h_S \har( |S|) \le C.
\end{equation}
\end{theo}
It follows from the upper bound in Matthews' theorem  and the
definition of the cover time that
\[
h_{\max} \le C \le h_{\max} \har(n).
\]
Thus, $h_{\max}$ approximates $C$ within a factor of $\har (n)
\approx \ln n$. Moreover, \corr{since $H(i,j)$ are} since the
$H(i,j)$'s are quite easily computable in polynomial time,
$h_{\max}$ is computable in polynomial time. Unfortunately,
$h_{\max}$ can be equal to the cover time (as shown by a path), as
well as a factor of $\ln n$ off the cover time (as shown by a
complete graph).

We could try to use the lower bound \raf{MATTHL} in Matthews'
Theorem, by maximizing over $S$. Unfortunately, this can be even
worse than
\corr{the upper bound.}
   $h_{\max}$.
For example, if $G$ consists of a
single edge with $N$ loops added at one of the nodes, then the
best Matthews lower bound is $1$, while the cover time, starting
from either the node with the loops or from the stationary
distribution, is about $N$. This problem is easy to fix: just
throw in the obvious lower bound $h_{\max}$. More precisely, let
$M_0$ be the maximum of $h_{\max}$ and the quantities $h_S \ln
|S|$ ($S\subseteq V$, $|S|\ge 2$). Then $M_0$ is a lower bound on
$C$, and (as we'll see) it is only a $(\ln\ln n)^2$ factor off. We
call $M_0$ the {\it augmented Matthews bound}.

A parameter closely related to the cover time is the {\it blanket
time}, introduced by Winkler and Zuckerman \cite{WZ}. The
definition provided below is a little bit stronger.
\corr{Omitted:  than the
definition of Winkler and Zuckerman.}

\smallskip

{\bf Definition.} Consider a random walk starting from a node $v$.
\corr{Let $r_T(x)$}
Let $r_{T,v} (x)$
be the number of visits to $x$
\corr{at time $T$.}
up to time $T$.
 Let ${\bf
B}$ be the first time $T$ when
\corr{the ratio ${r_T(i) /\pi_i \over
r_T(j)/\pi_j}\le 2$}
the ratio ${r_T(i) /\pi_i \over  r_T(j)/\pi_j}$ is at most 2
for any two nodes $i$ and $j$ (in particular,
all nodes are covered by this time). Let $B(v)$ be the expectation
${\bf B}$. The blanket time $B$ is the maximum of $B(v)$, over all
vertices $v$.

\smallskip

It is clear that $C \le B$.  Winkler and Zuckerman conjectured
that there is a constant $K$ so that $B \le K C$, and showed that
\[
B = O(C\ln n).
\]
The main goal of this paper is to improve the factor $O(\ln n)$ in
both problems mentioned above to $O((\ln \ln n)^2)$.

The following variant of the augmented Matthews bound $M_0$ is
\corr{the heart of}
 at the heart of
our study. Let  $\kappa(i,j)= H(i,j)+ H(j,i)$ be the commute time
between $i$ and $j$.
\corr{For any subset $S \subset V(G)$,}
For any $S \subset V(G)$,
let
$\kappa_S= \min_{i,j \in S} \kappa(i,j)$, and
\[
M = \max_{ S \subset V(G)} \kappa_S \ln |S|.
\]
As the following
\corr{proposition shows, the two bounds are essentially}
proposition shows, $M$ and $M_0$ are essentially
equivalent, but due to the symmetry of $\kappa$, $M$ will be
easier to handle.

\begin{prop}\label{TWOBOUNDS}~
\[
{1\over 8}M\le M_0 \le M.
\]
\end{prop}

Our main theorem is the following.

\begin{theo}\label{MAIN}
For every graph $G$ on $n$ vertices
\[
{1\over 2}M\le C \le B \le 10^5 M (\ln \ln n)^2
\]
\end{theo}
\corr{NEW}
(Of course the lower bound $C\geq M/8$ follows from
    Proposition 1.2 and Matthews' bound.)

It follows from this theorem that $M$ approximates $C$ within a
factor $O((\ln \ln n)^2)$ and $B \le K C (\ln \ln n)^2$, for some
constant $K$. It turns out, somewhat surprisingly, that both the
upper bound and the lower bound are sharp up to a constant factor.

The proof of the lower bound will give a somewhat stronger result.
Let $C(\pi)=\sum_i\pi_i C(i)$ denote the cover time when the walk
is started from a random node from
\corr{the stationary distribution.}
the stationary distribution $\pi$.

\begin{theo}\label{KAPPALOW}
For any graph $G$,
\[
C(\pi) \ge {1\over 2}M.
\]
\end{theo}

An important property of $M$ as an approximation of the cover time
is that it is efficiently approximable:

\begin{theo}\label{MAPPROX}
$M$ can be approximated within a factor of $2$ by a deterministic
polynomial algorithm.
\end{theo}

The rest of the paper is divided into five sections. In Section
\ref{APPROXM}, we describe an algorithm which computes $M$ up to a
factor of 2, proving Theorem \ref{MAPPROX}. In Section
\ref{FORMULAS},
\corr{as a preparation}
as preparation
 for the
 \corr{proof of the lower bound,}
 proof of the lower bound in Theorem~\ref{MAIN},
we derive some formulas for the cover time, which may be
interesting
\corr{on their own right.}
in their own right.
In Section \ref{LOWERBOUND}, we
complete the
\corr{proof of $M/2 \le C$,}
proof of Theorem~\ref{KAPPALOW},
 and also prove Proposition
\ref{TWOBOUNDS}. The proof of the
\corr{upper bound,}
upper bound in Theorem~\ref{MAIN},
which is the most
substantial part of this paper, follows in Section
\ref{UPPERBOUND}. In the final Section \ref{SHARPNESS}, we
\corr{give some  constructions}
give constructions
which show that
\corr{both the upper bound and lower bound}
      both the upper and lower bounds
\corr{in the Main Theorem}
in Theorem~\ref{MAIN}
can be attained.

\section{Approximating $M$}\label{APPROXM}

Since the commute times $\kappa(i,j)$ are polynomially computable,
the quantity $\kappa_S$ is also polynomially computable for any
set $S \subset V(G)$. However, the definition of $M$ involves all
(exponentially many) subsets of $V(G)$ and it is not clear that
one can compute $M$ in polynomial time. In the following, we show
that one can, at least, approximate $M$
\corr{within}
to within
a factor of 2 in
polynomial time.

A preliminary remark: the commute time $\kappa(i,j)$ satisfies the
triangle inequality:
\[
\kappa(i,k) \le \kappa(i,j) + \kappa(j,k),
\]
and hence we can consider it as a ``distance'' on the graph.

\smallskip

\noindent{\bf Algorithm.}  To start, pick  an arbitrary vertex
$v_1$. At the $i^{th}$ step ($i=1,2, \dots,n$), we have selected
the set $V_i = \{v_1, \dots, v_{i} \}$. Choose $v_{i+1}$ to be a
vertex $v\in V\setminus V_i$ whose distance $\min_{u \in V_i}
\kappa(u,v)$ from $V_i$ is maximum. Compute $M_i = \kappa_{V_i}
\ln i$ for all $i=2,3, \dots, n$, and output $M' = \max_{i} M_i$.

\smallskip

Since $\kappa(i,j)$ are polynomially computable, our algorithm
runs in polynomial time. Moreover, $M' \le M$ by definition. It
remains to show that $2M' \ge M$.
\corr{ Omitted:, independently from the choice of
$v_1$.}

Assume that $M$ is attained at a set $S \subset V(G)$ of
cardinality $s$. We claim that $2M_s \ge M$. It suffices to show
that $\kappa_{V_s} \ge \kappa_{S}/2$.

Let $R = S \setminus V_s$.  If $R$ is empty then $S = V_s$ and we
are done, so we assume that $|R|=r >0$. By the description of the
algorithm, $\kappa_{V_s} = \kappa(v_s, v_j)$ for some $j < s$.

For each vertex $x \in R$, there is a vertex $y_x \in V_{s-1}$ so
that $\kappa(x,y_x) \le \kappa_{V_s}$. If $y_x\in S$ for some $x$,
then $\kappa_S\le \kappa(x,y_x) \le \kappa_{V_s}$, and we are
done. If $y_x\in V_{s-1}\setminus S$ for all $x\in R$, then (using
that $|V_{s-1}\setminus S|=(s-1)-(s-r)=r-1<r$) the pigeon hole
principle gives that there are $x$ and $x'$ in $R$ so that $y_x =
y_{x'}=y$. So by the triangle inequality
\[
\kappa(x,x') \le \kappa(x,y) + \kappa(y,x') \le 2 \kappa_{V_s}.
\]
By definition $\kappa(x,x') \ge \kappa_S$ and the proof is
complete.

\smallskip

\noindent{\bf Remark.} The only property of the commute times we
use here is the triangle inequality.  Therefore, our result holds
in a more general setting. Consider a metric $w$ on a finite set
$V$ of $n$ points. For any subset $S\subset V$, let $w_S =
\min_{i,j \in S} w(i,j)$ (if $S$ has less than 2 elements,
$w_S=0$). Define
\[
W = \max_{S \subset V} w_S f(|S|),
\]
where $f$ is any non-negative function defined on the set of
non-negative integers.

\begin{cor}\label{WEIGHTS}
\corr{For any finite metric space,} For any finite metric space
and any non-negative $f$, the above algorithm (with $\kappa(i,j)$
replaced by $w_{ij}$) computes $W$ within a factor of 2.
\end{cor}

\section{Formulas for the cover time}\label{FORMULAS}

Fix a set $S\subseteq V$, $|S|=s \ge 2$, and a starting node $v$.
For a given random walk $(v=v^0,v^1, v^2, \dots)$, and a set
$T\subseteq S$, let $Z(T)$ denote the set of nodes of $S$ not seen
before $T$ is first reached. Thus $T\subseteq Z(T)$. Define, for
$i,j\in S$,
\[
A(i,j)=\cases{{1\over |Z(i)|(|Z(i)|-1)},& if $j\in Z(i)\setminus
\{i\} $,\cr
                0,& otherwise,\cr}
\]
(this number depends on the walk) and let $a(i,j)=\E[A(i,j)]$. We
have
\[
\sum_{i\in S}\sum_{j\in S} A(i,j)=\sum_{i\in S} \sum_{j\in
Z(i)\setminus \{i\}} {1\over |Z(i)|(|Z(i)|-1)}
\]
\[
=\sum_{i:\,|Z(i)|>1} {1\over |Z(i)|} = {1\over 2}+{1\over
3}+\dots+{1\over s},
\]
and thus
\begin{equation}\label{ASUM}
\sum_{i\in S}\sum_{j\in S} a(i,j)= {1\over 2}+{1\over
3}+\dots+{1\over s}\approx \ln s.
\end{equation}

Using this notation, we can state a formula for the expected
number $C(v,S)$ of steps
\corr{before}
until all nodes of $S$ are
\corr{visited:}
visited. The basic idea here is
similar to that in Mathhews' theorem.

\begin{lemma}\label{HITTING}~
\[
C(v,S)= {1\over s}\sum_{j\in S} H(v,j) + \sum_{i\in S}\sum_{j\in
S} H(i,j)a(i,j).
\]
\end{lemma}

\Proof Let $X_k$ be the number of steps
\corr{before visiting the $k$-th
node in $S$.}
required to see $k$ nodes of $S$.
 Clearly $C(v,S)=\E[X_s]$. Let $T(i)$ be the number of
steps
\corr{before first visiting}
required to see
node $i$.
The following algebraic
identity is easy to verify:
\begin{eqnarray}\label{AAA}
X_s &=& {1\over s} \sum_{k=1}^s X_k \nonumber
\\ &+& \sum_{1\le k < m\le s}
{1\over (s-k)(s-k+1)} (X_m-X_k)
\end{eqnarray}
Now here
\[
\E\Bigl[\sum_{k=1}^s X_k\Bigr]=\sum_{i\in S} \E[T(i)] = \sum_{i\in
S} H(v,i).
\]
\corr{For the second sum,}
For the second sum in (\ref{AAA}),
we fix the first $X_k$ steps, then
\[
\sum_{m:\,k \le m\le s} (X_m-X_k) = \sum_{j\in Z(v^{X_k}) }
(T(j)-T(v^{X_k})),
\]
and hence
\[
\E\Bigl[\sum_{m:\,k \le m\le s} (X_m-X_k)\Bigr] = \sum_{j\in
Z(v^{X_k})} H(v^{X_k},j).
\]
Summing over $k$, we get
\[
\sum_{k=1}^{s-1} {1\over (s-k)(s-k+1)}\sum_{j\in Z(v^{X_k})}
H(v^{X_k},j)
\]
\[
= \sum_{i\in S} {1\over (|Z(i)|-1)|Z(i)|}\sum_{j\in Z(i)} H(i,j)
\]
\[
= \sum_{i\in S}\sum_{j\in S} H(i,j) A(i,j).
\]
Taking expectation again, we get the lemma.\proofend

For $i,j\in S$, let
\[
Q(i,j)={1\over |Z(ij)|(|Z(ij)|-1)},
\]
and $q(i,j)=a(i,j)+a(j,i)=\E[Q(i,j)]$. Using the identity
\begin{equation}\label{HITDIFF}
H(\pi,j) - H(\pi,i) = H(i,j)- H(j,i).
\end{equation}
due to Tetali and Winkler \cite{TW}, which implies that
\begin{equation}\label{HITKAPPA}
H(i,j)= {1\over 2} \kappa(i,j) + {1\over 2}(H(\pi,j) - H(\pi,i)),
\end{equation}
a simple computation gives the following lemma:

\begin{lemma}\label{KAPPA}~
\begin{eqnarray}\label{CS}
C(v,S)&=& {1\over s}\sum_{j\in S} (H(v,j)-H(\pi,j))\\ &+& {1\over
4} \sum_{i\in S}\sum_{j\in S} \kappa(i,j)q(i,j).
\end{eqnarray}
\end{lemma}

Let $q_{\pi}(i,j)$ be the expectation of $q(i,j)$, when the
starting node $v$ is chosen at random from the stationary
distribution. Averaging over $v$, the first term in \raf{CS}
cancels, and we get

\begin{cor}\label{FROMPI}~
\[
C(\pi,S)={1\over 4} \sum_{i\in S}\sum_{j\in S}
\kappa(i,j)q_{\pi}(i,j).
\]
\end{cor}

\section{Proof of the lower bound.}\label{LOWERBOUND}

\noindent{\bf Proof of Theorem \ref{KAPPALOW}} This
\corr{follows from
Corollary \ref{FROMPI} and (\ref{ASUM}) easily:}
follows easily from Corollary \ref{FROMPI} and (\ref{ASUM}):
\[
C(\pi)\ge C(\pi,S)={1\over 4} \sum_{i\in S}\sum_{j\in S}
\kappa(i,j)q_\pi(i,j)
\]
\[
\ge {1\over 4} \kappa_S \sum_{i\in S}\sum_{j\in S} q_\pi(i,j)=
{1\over 2} \kappa_S \ln |S|.
\]

\noindent{\bf Proof of Proposition \ref{TWOBOUNDS}}. It is obvious
that for every $S\subseteq V$
\[
 h_S \ln |S|\le  {1\over 2}\min_{i,j\in S, i\not=j}\kappa(i,j)
  \ln |S| \le {1\over 2}M,
\]
and for every $i,j\in V$,
\[
H(i,j)\le \kappa(i,j) \le M.
\]
Hence $M_0\le M$.

To prove the other bound, let $S$ be the set attaining the maximum
in the definition of $M$. If $h_{\max} >M/4$, then we have nothing
to prove, so suppose that $H(i,j) \le M/4$ for all $i$ and $j$.

We define a digraph $D$ on $S$ as follows. There is an edge from
$i$ to $j$ if and only if $H(i,j) \le \kappa(i,j)/4$. In this
case, it is clear that $H(j,i)-H(i,j) \ge  \kappa(i,j)/2 \ge
\kappa_S/2$.

Let $i_0i_1 \dots i_m$ be a (directed) path of length $m$ in $D$.
Then by the cycle law \cite{TW} we have
\[
H(i_m, i_0) \ge \sum_{l=0}^{m-1} (H(i_{l+1},i_l) - H(i_l,i_{l+1}))
\ge m \kappa_S/2.
\]
On the other hand, $H(i,j) \le \kappa_S \ln |S|/4$ for all $i,j
\in S$. This implies that $m \le \ln |S|/2$. A theorem of Gallai
\cite{Gal} implies that $D$ is $\ln |S|/2$ colorable, and
therefore $D$ contains an independent set $I$ of size at least
$2|S| /\ln |S|>|S|^{1/2}$. By the definition of $D$,
\[
H(i,j) \ge  \kappa(i,j)/4 \ge \kappa_S/4
\]
for any $i,j \in I$. Therefore,
\[
\min_{i,j \in I} H(i,j) \ln |I|  \ge {1\over 4}\kappa_S  {1\over
2}\ln |S| = {1\over 8} M.
\]
\proofend

\section{Proof of the upper bound}\label{UPPERBOUND}

We need a Chernoff type large deviation inequality, which will be
shown using fairly standard arguments.

\begin{lemma}
\label{geo} Let $X_1, ..., X_k$ be independent non-negative
integer valued random variables with $$ \Pr [X_i = m] \le a
(1-p)^m \quad \forall \, m \geq 1 $$ for some numbers $a>0$ and
$0< p< 1$. Let $X=\sum_{i=1}^{k} X_i$. Then for
 any $L>0$ $$ \Pr [ X  - \E[X ]  \le - L ]
\le \exp\Big( - \frac{p^3 L^2}{4(1-p)ak} \Big). $$
\end{lemma}

\Proof As usual, we first estimate
$\E[e^{-\lambda  X_i}]$ for $\lambda >0$. Taylor expansion
gives
\begin{eqnarray*}
\E[e^{\lambda  X_i} ] &=& 1-\lambda  \E[X_i] + (\lambda ^2/2)
\E[X_{i}^{2} e^{-\lambda ^* X_i}]\\ &\le& \exp(-\lambda  \E[X_i] +
(\lambda ^2/2) \E[X_{i}^{2} ])
\end{eqnarray*}
for some $\lambda^*$ between $0$ and $\lambda$.
Since
\begin{eqnarray*}
\E[X_{i}^{2} ] &=& \sum_{m=1}^{\infty} m^2
\Pr[ X_i = m ]\\ &\le& a \sum_{m=1}^{\infty} m^2
(1-p)^m
\end{eqnarray*}
and
$$\sum_{m=1}^{\infty} m^2 (1-p)^m = \frac{(1-p) + (1-p)^2
}{p^3} \leq
\frac{2(1-p)}{p^3},$$ it follows that
 $$ \E[e^{-\lambda  X_i} ]\le \exp \Big(-\lambda  \E[X_i] +
\frac{a(1-p)\lambda ^2}{p^3} \Big). $$
\noindent Therefore
\begin{eqnarray*}
\E[ e^{-\lambda  (X - \E[X] )}] &=& \prod_i \E[e^{-\lambda  (X_i
-\E[X_i])} ] \\ &\le& \exp \Big( \frac{ak(1-p)\lambda
^2}{p^3}\Big).
\end{eqnarray*}
\noindent Taking $\lambda =
 p^3 L /(2(1-p)ak)$, we have
\begin{eqnarray*}
\Pr [ X-\E[X] \le  - L ] &\le& \E[ e^{-\lambda (X-E[X]+L)}] \\
&\leq& \exp\Big(-\lambda L +\frac{ak(1-p)\lambda
^2}{p^3}\Big) \\
&=& \exp\Big(-\frac{p^3 L^2  }{4(1-p)ak} \Big) .
\end{eqnarray*}
\proofend

\begin{lemma}\label{BACKFORTH}
Let $i$ and $j$ be two nodes and $k\ge 1$. Let $W_k$ be the number
of times $j$ had been visited when $i$ was visited the $k$-th
time. Then for every $\eps>0$,
\[
\Pr\Bigl[W_k<(1-\eps){\pi_j\over\pi_i}k\Bigr] \le
\exp\left(-\eps^2k\over 4\pi_i\kappa(i,j)\right)
\]
\end{lemma}

\Proof Let us restrict the Markov chain to $i$ and $j$ only. It is
well known that we get a time-reversible Markov chain with
transition probabilities
\[
\matrix{\hat p_{ii}=1-{1\over \pi_i\kappa(i,j)}&
\hat{p}_{ij}={1\over \pi_i\kappa(i,j)}\cr \hat{p}_{ji}={1\over
\pi_j\kappa(i,j)}& \hat{p}_{jj}=1-{1\over \pi_j\kappa(i,j)}\cr }
\]
and stationary probabilities
\[
\hat\pi_i = {\pi_i\over \pi_i+\pi_j}, \qquad \hat\pi_j =
{\pi_j\over \pi_i+\pi_j}.
\]
We may consider this very simple Markov chain to prove the lemma.

Define $X_k$ to be the number of visits to $j$ during the $k^{\rm
th}$ return trip from $i$ to itself, that is,
\[
X_j = W_{k+1} - W_k.
\]
It is clear that the $X_j$ are i.i.d. with
\begin{eqnarray*}
\Pr [ X_1 = 0] &=& \hat{p}_{ii}
\\ \Pr [ X_1 = m] &=& \hat{p}_{ij}\hat{p}_{jj}^{m-1} \hat{p}_{ji}
\\ \E[X_1]&=& {\pi_j\over\pi_i}
\end{eqnarray*}
Clearly $W_k=X_1 + \cdots + X_k$. Thus
$$
\E[W_k] = k \E[X_1] = k {\pi_j \over \pi_i}.
$$
\newcommand{\pp}{\hat{p}}
Applying Lemma \ref{geo} with $a=\hat{p}_{ij}\hat{p}_{ji}
(1-\hat{p}_{ji})^{-1}$, $p= \hat{p}_{ji}$ and
$L=\eps{\pi_j\over\pi_i}k$, we obtain
\[
\Pr [W_k<(1-\eps){\pi_j\over\pi_i}k] =\Pr \Big[X-\E[X] \le -
\eps{\pi_j\over \pi_i}k\Big]
\]
\[
\le \exp\Big( - \frac{\pp_{ji}^3
\eps^2\pi_i^2k^2}{4\pi_j^2\pp_{ij}\pp_{ji}k}\Big) = \exp\Big( -
\frac{\eps^2k}{4\pi_i\kappa(i,j)}\Big)
\]
\proofend

\noindent{\bf Proof of the upper bound in Theorem \ref{MAIN}}.
Consider the ordering $(v_0,v_1,\dots,v_{n-1})$ of the nodes of
$G$ as obtained by the Algorithm in section \ref{APPROXM}. For
convenience, relabel the nodes by $(1,\dots,n)$. Recall that each
$i>1$ is a node farthest away from the set $\{1,\dots,i-1\}$ in
distance $\kappa$.

For each node $i>1$, let $i'$ be a node with $i'\le\sqrt{i}$ and
$\kappa(i,i')$ minimal. Clearly, the edges $ii'$ form a tree
$\TT$. We consider $1$ as the root of the tree. It is also clear
that the depth $d$ of $\TT$ is at most $1.5 \ln\ln n$.

Our next observation is that
\begin{equation}\label{TREEKAPPA}
\kappa(i,i')\le {2M\over\ln i}.
\end{equation}
Indeed, let $S=\{1,\dots,\lfloor\sqrt{i}\rfloor+1\}$. Then, by the
definition of $M$,
\[
\kappa_S\le {M\over \ln |S|} < {2M\over\ln i},
\]
and hence there exist nodes $u,v\in S$ with $\kappa(u,v)< 2M/\ln
i$. We may assume that $u\le v$. By the choice of the ordering,
there exists a node $j\le v-1 \le \sqrt{i}$ such that
$\kappa(i,j)\le\kappa(u,v)$. It follows that
\[
\kappa(i,i')\le\kappa(u,v)\le  {2M\over\ln i}.
\]

Set $\eps=1/(8\ln\ln n)$ and $T_0= \lceil 400M/\eps^2\rceil$. Our
next goal is to bound the probability that ${\bf B}>T$ for some
$T\ge T_0$.

Set $F(i)=r_T(i)/(T\pi_i)$. On the average, $F(i)=1$. If the event
``${\bf B} > T$" occurs, then there exists an edge $ii'$ of $\TT$
with one of the following properties:

\smallskip

(A) $F(i') \geq 0.9(1+\ln i')^{-\eps}$ and $F(i)< 0.9(1+\ln
i)^{-\eps}$;

\smallskip

(B) $F(i') \leq 1.1(1+\ln i')^{\eps}$ and $F(i) > 1.1(1+\ln
i)^{\eps}$;

\smallskip

(C) $F(i') \leq 0.9(1+\ln i')^{\eps}$ and $F(i)> 0.9(1+\ln
i)^{\eps}$;

\smallskip

(D) $F(i') \geq 1.1(1+\ln i')^{\eps}$ and $F(i)< 1.1(1+\ln
i)^{\eps}$.

\smallskip

Indeed, if $\bf{B} $ is larger than $T$, then there exists a node
$u$ such that either $F(u)>\sqrt{2}$ or $F(u)<1/\sqrt{2}$. Suppose
that e.g. the second occurs. We assume that $n>10$, to exclude
some trivial complications. Then $F(u)< 0.9(1+\ln u)^{-\eps}$. We
also know that there is a node $v$ with $F(v)>1>0.9(1+\ln
w)^{\eps}$. If $F(1) > 0.9$, then along the path from $u$ to $1$
there is an edge with property (A). If  $F(1) \le 0.9$, then along
the path from $v$ to $1$ there is an edge with property (C).

We call such an edge ``bad". To bound the probability that an edge
is bad, we have to bound the probabilities of (A), (B), (C) and
(D) separately. This is very similar in all cases, and we give the
details for (A). Let $k=\lceil 0.9(1+\ln i')^{-\eps}\pi_{i'}
T\rceil$, and consider the step when $i'$ is reached the $k$-th
time. By (A), the number of times we have seen $i$ is
\[
W_k<0.9(1+\ln i)^{-\eps} \pi_i T \le \left({1+\ln i\over 1+\ln
i'}\right)^{-\eps} {\pi_i\over \pi_{i'}}\cdot k
\]
\[
< 2^{-\eps}{\pi_i\over \pi_{i'}}\cdot k < \left(1-{\eps\over
4}\right) {\pi_i\over \pi_{i'}}\cdot k ,
\]
and hence by Lemma \ref{BACKFORTH},
\[
\Pr[A] \le \exp\left(-\eps^2 k\over
100\pi_{i'}\kappa(i,i')\right).
\]
Now here
\[
{k\over\pi_{i'}} \ge 0.9(1+\ln i')^{-\eps}T \ge {1\over 2}T
\]
and hence by \raf{TREEKAPPA},
\[
\Pr[A]<\exp\left(-\eps^2 T \ln i \over 200M\right)<
i^{-T\eps^2/(200M)}.
\]
The probability that this happens for some edge $ii'$ is at most
\[
\sum_{i=2}^n i^{-T\eps^2/(200M)}< 2\cdot 2^{-T\eps^2/(200M)},
\]
(using here that $T\ge T_0$) and hence the probability that ${\bf
B}>T$ is at most $8\cdot 2^{-T\eps^2/(200M)}$. Thus
\begin{eqnarray*}
\E[{\bf B}] &=& \sum_{T=0}^{\infty} \Pr[{\bf B}> T] \le T_0 +
\sum_{T=T_0}^{\infty} \Pr[{\bf B}> T]
\\ &\le& T_0+ 8\sum_{T=T_0}^{\infty} 2^{-T\eps^2/(200M)} < 2T_0,
\end{eqnarray*}
which proves the theorem.
\proofend

\vskip 1mm

{\bf \noindent Remark.} We may prove the upper bound using a slightly
 different approach. Let  $S_0$ be the set of all vertices
 and inductively define $S_i$ to be a maximal subset
of $S_{i-1}$ such that $\kappa_{S_i}  >  \kappa_i := 2^i M /\log
n$. If such a subset does not exist, $S_i$ consists of a vertex in
$S_{i-1}$ and the construction stops. Since $\kappa_{S_i}  \leq M$
unless $|S_i |=1$, this procedure stops within $ O(\ln \ln n) $
steps. The advantage of this approach is that we may have a better
upper bound if the procedure stops earlier. For example, if $G$ is
a complete graph, then the construction stops after $1$ step. More
generally, for each $x \in S_{i}\setminus S_{i+1}$, take a vertex
$y\in S_{i+1}$ with $ \kappa(x,y) < \kappa_{i+1}$. This is
possible since $S_{i+1}$ is a maximal subset. Regarding the pair
$xy$ as an edge, this gives a tree with depth at most $O(\ln \ln
n)$. Let $l$ be the minimum possible depth. Then the same proof
would yield $B = O(l^2 M)$.

\section{The sharpness of Main Theorem}\label{SHARPNESS}

In this section we show that both the lower bound and upper bound
in the Main Theorem \ref{MAIN} are sharp, up to a constant factor.
More exactly, we give an example where $B$ and $C$ are of order
$\Theta (M)$  and also one where $B$ and $C$ are of order $\Theta(
M(\ln\ln n)^2)$.

The proof for the lower bound is easy: for the complete graph on
$n$ vertices, all three parameters $B, C$ and $M$ are
\corr{of order}
$\Theta (n \ln n)$.

The construction to match to upper bound is more complicated. It
is a tree of depth $d$ defined as follows. The root is at level
$1$. Each vertex at the $i^{th}$ level has $2^{2^i}$ children, and
the edge between the mother and a child has multiplicity $2^i$.

The number of vertices in the $i^{th}$ level is
\[
N_i = \prod_{j=1}^{i-1}  2^{2^{j}} = 2^{ 2^i -2}.
\]
The number of the vertices in the whole tree is
\[
n=\sum_{i=1}^d N_i =  \sum_{i=1}^d 2^{2^i-2} = N_d(1+o(1)).
\]
The number of edges between the $i^{th}$ and $(i+1)^{th}$ level is
$E_i= 2^i N_{i+1} = 2^{ 2^i + i -2}$. The total number of edges is
\[
E= \sum_{i=1}^{d-1} E_i = E_{d-1} ( 1+ o(1)).
\]
Notice that $d = \Omega(\ln\ln n) $.
We first show
$$ M= \Theta(E). $$
It is well-known that the commute time between two vertices
$x,y$ in a tree (possibly with multiple edge) is
\begin{equation}\label{treecom}
 2E\cdot  \sum_{j=0}^{l-1} \frac{1}{m(x_j x_{j+1})} ,
 \end{equation}
where $x=x_{_0} ,..., x_{_l}=y$ is the path connecting $x$ and
$y$,  and $m(vw)$ is the multiplicity of the edge $vw$.
For a lower bound, consider the set $S_2$ of (four)
vertices in level $2$. The commute time is $2E$
for any pair (by (\ref{treecom}),
which gives $ M \geq 2E \ln 4$.
For an upper bound, let $S$ be a set of size at least $2$.
Then take the maximum level $i_{_0}$ such that there is a vertex
in level $i_{_0}$ having  at least
two descendants (including itself) in $S$.
Since the multiplicity of an edge
 geometrically increases as the level increases,
 (\ref{treecom}) implies that
the commute time of
a pair who has a common ancestor in level $i_{_0}$   is at most
$O(E/2^{i_{_0}} )$, especially $ \kappa_{S} = O(E/2^{i_{_0}} )$.
Moreover, since no pair has a common ancestor in level $i_{_0}+1$,
the number of vertices in $S$ below level $i_{_0}$ is at most
$N_{i_{_{0}}+1}$. Trivially, the number vertices  of $S$ above or
in level $i_{_0}$ is at most $\sum_{i=1}^{i=i_{_0}} N_i  =
o(N_{i_{_{0}}+1})$. Thus $|S| = (1+o(1))N_{i_{_{0}}+1}$ and
$$ \kappa_S \ln |S| = O(E). $$


In the rest of this section, we  shall omit
\corr{Omitted: unnecessary}
floors and
ceilings, for the sake of a clearer presentation.

\begin{claim}\label{TREEC}
The cover time of this tree satisfies $C = \Omega (M d^2)$.
\end{claim}

\noindent{\bf Proof.} It suffices to show that for a sufficiently
large constant $K$, a walk of length $T= d^2E/K$, starting from a
stationary point, covers the tree with probability at most $1/2$.

 To
start, set $k= 10\ln \ln d$ and define a sequence $b_i$ as follows
\[b_k = d^2/\sqrt{K}, b_{i} =  b_{i-1} ( 1- \sqrt {\frac{1}{2} / b_{i-1}}),
\]

\noindent for all $i >k$. Let $l$ be the first index such that
$b_l \le 1/2$. Arithmetic shows that if $K$ is sufficiently large
then $l < d-1$. Set $a_i = 2^i b_i$ and $m_i = 2^{2^{i+1}}$, a
simple calculation shows

\[a_i = 2 a_{i-1} - \sqrt { \frac{1}{4} a_{i-1} \ln m_{i-1}} . \]

Let  $X_i$ denote the minimum number of times a multi-edge from
level $i$ to level $i+1$ is crossed in a finite walk. We say that
a walk is a $T_i$-walk if it stops when $X_i= a_i$ and denote by
$A_i$ the event that a $T_i$ walk  covers the tree. Furthermore,
let  $B$  be the event that  a walk of length $T$ satisfies $X_k
\ge a_k$. Notice that

\[ \Pr (\hbox{ A walk of length $T$ covers the tree}) \]
\[\le \Pr(\hbox{B} ) + \Pr (A_k) \]

 The expectation of the number of crosses of any multi-edge between the
 $k^{th}$ level and the $(k+1)^{th}$ level
  is $2^k T/E= 2^kd^2/K$, where $2^k$ is th multiplicity of the edge.
   On the other hand, $a_k = 2^k d^2/\sqrt{K}$ by definition.
 Therefore, by Markov's inequality $\Pr(B)$ is at most
  $1/\sqrt{K} < 1/3$. To finish the proof, we show that
 $\Pr (A_k) =o(1)$. Observe that for any $i \ge k$, $\Pr(A_i)$ is
 upper bounded by

\[
 \Pr (\hbox{ a $T_i$-walk satisfies $X_{i+1} \ge a_{i+1}$} )
  + \Pr ( A_{i+1} ). \]

\noindent It follows that

\[ \Pr (A_k) \le
\sum_{i=k}^{l-1} \Pr (\hbox{ a $T_i$-walk satisfies $X_{i+1} \ge
a_{i+1}$} ) \] \[+  \Pr ( \hbox{a $T_l$-walk covers the tree} ).
\]

\noindent To show that $\Pr (A_k)= o(1)$, it now suffices to prove
that

\begin{equation}\label{TIWALK}
\sum_{i=k}^{l-1}  \Pr (\hbox{ a $T_i$-walk satisfies $X_{i+1} \ge
a_{i+1}$} ) = o(1) \end{equation} and

\begin{equation}\label{TLWALK}  \Pr ( \hbox{a $T_l$-walk covers the tree}
)=o(1).
\end{equation}

It will be useful to think about the walk using a ``balls and
urns" model. Consider a vertex $u$ on level $i$. Attach to each
neighbor of $u$. Any time we exit node $u$, drop a ball into the
corresponding urn. Then balls will be dropped into the urns
independently, so that the urns corresponding to the children of
$u$ have the same probability, and the urn corresponding the
parent of $u$ has half this probability. Conversely, if for each
node, we decide how to drop balls into the urns, then we determine
a unique walk. It is important to notice that the number of times
an edge is crossed depends only on the ball distributions
corresponding to nodes above the edge.

Assume that the multi-edge between $u$ and its parent $v$ is
crossed $x$ times; then the numbers of crossing of the multi-edges
going down from $u$ is the same as the number of balls in the big
urns at the moment the small urn has $x$ balls.

Using the balls and urns terminology,  \raf{TIWALK} follows from
the following lemma.

\begin{lemma}\label{URN}   Assume that  $a$ and $m$ are large numbers,
and $a'= 2a - \sqrt{ \frac{1}{4} a \ln m} \ge 0$. Drop balls into
one small urn and $m$ big urns until the small urn has $a$ balls,
then with probability at least $1- \Big(\exp (-\ln^{2/3} m )
+\exp(-m^{1/2 }) \Big) $, one of the big urns has at most $a' $
balls.
\end{lemma}

{\noindent \bf Proof.} We use the following fact which is easy to
prove. If $X$ is sum of i.i.d. binary random variables and $X$ has
large expectation $\mu$, then for any $\sqrt{\mu} \le  L \le  \mu$

\begin{equation}\label{CONCEN}  \exp( - 2L^2/ \mu) \le
\Pr ( X \le \mu - L)  \le 2\exp (- L^2/ 2\mu).
\end{equation}

To prove the lemma,  we first show that with  probability at least
$ 1- \exp (\ln^{2/3} m) $, at the first moment when the small urn
has $a$ balls, the number of balls dropped is at most $A=a(2m+1) +
4m\sqrt {a \ln^{2/3} m} $. To show this, it is enough to prove
that if one drops $A$ balls randomly into one small urn and $m$
big urns, then with probability at least $1- \exp(-\ln^{2/3} m)$
the small urn has at least $a$ balls. The number of balls in the
small urn can be expressed as a sum of $A$ i.i.d. binary random
variables and has expectation $\mu= A/(2m + 1) = a +L$, where $L=
2\sqrt {a \ln^{2/3} m} +o(1)$. The claim follows directly from the
upper bound in \raf{CONCEN}, with room to spare.

To finish the proof of the lemma, we show that if we drop $A$
balls into one small urn and $m$ big urns, then there is a big urn
with at most $a'$ balls with probability at least  $1-
\exp(-m^{1/2+ o(1)})$. The number of balls in a fixed big urn is a
sum of $A$ i.i.d. binary random variables and has  expectation
$A/(m+1/2) $. Set  $L'= A/ (m +1/2) - a'$; it is clear that $L' =
(1+o(1)) \sqrt{ \frac{1}{4} a \ln m}$. We say an urn is ``good" if
it has at most $a'$ balls and ``bad" otherwise. By the lower bound
in \raf{CONCEN}, the probability that a fixed urn is ``good" is at
least

\[p =\exp \Big(- 2 L'^2/ (A/(m+1/2)) \Big)  \ge m^{-1/2 }.\]

So the probability that an urn is ``bad" is at most $1-p$. Observe
that the events ``urn $U_1$ is bad" and ``urn $U_2$ is bad" are
negatively correlated, for any two fixed urns $U_1$ and $U_2$.
Using FKG inequality and induction, we can show

\[ \Pr (\hbox{ all $m$ urns are ``bad"}) \le (1-p)^{m} \]
\[ \le (1- m^{1/2 }) ^{m} \le \exp (-m^{1/2 }), \]

\noindent concluding the proof. \proofend

Now \raf{TIWALK} follows from the previous lemma and the fact that

\[\sum_{i=k}^{l-1} \exp (-\ln^{2/3} m_i )+ \exp (-m_i^{1/2 +o(1)}) = o(1). \]

\noindent Here we need to use the condition  $k = 10\ln \ln d$.

To prove \raf {TLWALK}, it suffices to prove show that
 if one drops balls into one small urn and $m_l= 2^{2^{l+1}}$ big urns
 until the small urn has
$a_l \le 2^l/2$ balls, then with probability at least $1-o(1)$,
there is an empty big urn. Similar to the proof of Lemma
\ref{URN}, one can show that at the time when the small urn has
$a_l$ balls, with probability $1-o(1)$, at most $3 a_l m_l$ balls
have been dropped (the constant $3$ is generous). To conclude,  we
show that if we drop $A_l= 3a_l m_l$ balls into $m_l$ identical
urns, then with probability $1-o(1)$, there is an empty urn. Since
$a_l \le 2^l/2$ and $m_l=2^{2^{l+1}}$, $A_l \le \frac{2}{3}m_l \ln
m_l $, the claim  follows  by a standard coupon collector
argument. \proofend


\begin{thebibliography}{99}

\bibitem{AF}
D.~J.~Aldous and J.~Fill, {\it Time-reversible Markov chains and
random walks on graphs} (book in preparation).



\bibitem{Gal}
T.~Gallai, On directed paths and circuits, in: {\it Theory of
Graphs} (Proc. Colloq., Tihany, 1966) Academic Press, New York
(1968), 115--118.

\bibitem{Mat}
P.~Matthews, Covering problems for Brownian motion on spheres,
{\it Ann. Prob.} {\bf 16} (1988), 189--199.

\bibitem{TW}
P.~Tetali and P.~Winkler, Simultaneous reversible Markov chains,
in: {\it Combinatorics, Paul Erd\H os is Eighty}, Vol. 1 (ed.
D.~Mikl\'os, V.~T.~S\'os, T.~Sz\H onyi), Bolyai Society, Budapest,
1993, 422-452.

\bibitem{WZ}
P.~Winkler and D.~Zuckerman, Multiple cover time, {\it Random
Structures Algor.} {\bf 9} (1996), 403--411.

\end{thebibliography}
\end{document}